\documentclass[12pt]{article}

\newcommand{\C}{\mathbf{C}}

\begin{document}

\title{Multidimensional analogues of Bohr's theorem on power
  series\thanks{Middle East Technical University Preprint
    Mathematics 164/98. This work was supported by the BSF
    grant No 94-00113.}}

\author{Lev Aizenberg}

\date{March 30, 1998}

\maketitle

\newtheorem{teo}{Theorem}
\begin{abstract}
  Generalizing the classical result of Bohr, we show that if an
  $n$-variable power series converges in an $n$-circular bounded
  complete domain~$D$ and its sum has modulus less than~$1$, then
  the sum of the maximum of the moduli of the terms is less
  than~$1$ in the homothetic domain $r \cdot D$, where $r = 1-
  \sqrt[n]{2/3}$.  This constant is near to the best one for the
  domain $ D = \{z: |z_1 |+ \ldots + |z_n | < 1 \}$.
\end{abstract}

\section{Preliminaries}
It is known the following classical Bohr's result \cite{B}, which
final form was contributed by M.~Riesz, I.~Schur and N.~Wiener.
\begin{teo}
If a power series
$$\sum_{k=0}^\infty c_k z^k \eqno{(1)}$$ converges in the unit
disk and its sum has modulus less than $1$ then
$$\sum_{k=0}^\infty | c_k z^k | < 1 \eqno{(2)}$$
in the disk $\{z:  |z| < 1/3 \}$ and 
the constant $1/3$ cannot be improved.
\end{teo}

Recently H.~P. Boas and D.~Khavinson obtained some
multidimensional generalization of this result (\cite{BK}, see
also for more references there).  Denote by $K_n$ the largest
number such that if the series
$$\sum_\alpha c_\alpha z^\alpha \eqno{(3)}$$
converges in the unit polydisk 
$U_1 = \{z:  |z_j| < 1$,  $j=1, \ldots, n \}$ 
and the estimate
$$\biggl| \sum_\alpha c_\alpha z^\alpha \biggr| < 1,\eqno{(4)}$$
is valid there, then 
$$\sum_\alpha | c_\alpha z^\alpha | < 1, \eqno{(5)}$$
holds in $K_n \cdot U_1$; here  $\alpha = ({\alpha}_1,  \ldots, 
{\alpha}_n)$, all ${\alpha}_j$ 
are non-negative integers,
$z = (z_1, \ldots, z_n)$, 
$z^{\alpha} = z_1^{{\alpha}_1} \cdots
z_n^{{\alpha}_n}$.
\begin{teo} [Boas, Khavinson]
It is true for $n>1$ that
$$\frac{1}{3 \sqrt{n}} < K_n < \frac{2 \sqrt{\log n}}{\sqrt{n}}
\eqno{(6)}$$
\end{teo}
\begin{teo} [Boas, Khavinson] 
Let the series (3) converge in a
complete $n$-circular domain (Reinhardt's domain) $D$ and 
(4) holds in $D$. Then (5) is true in the homothetic domain
$K_n \cdot D$.
\end{teo}

Notice that Remark 1 from \cite{BK}, in fact, contains a result
stronger then the left part of the inequality (6), namely
$$K_n > \left\{\begin{array}{lll}\frac{2}{5 \sqrt{n}} 
& \mbox{for} & n > 1 ,\\& & \\
\frac{1}{2 \sqrt{n}} & \mbox{for} & \mbox{large enough}  \; n.
\end{array} \right.$$

\section{Main Results}
We consider some other multidimensional variants of Bohr's problem.
Denote by $B_n(D)$ the biggest number $r$ such that if the 
series (3) converges in a complete $n$-circular domain $D$ 
and (4) holds in it then
$$\sum_\alpha \sup_{D_r} | c_\alpha 
z^\alpha | < 1, \eqno{(7)}$$
where $D_r = r\cdot D$ is the homothetic transformation of $D$. 
If $D = U_1$ then $B_n(D) = K_n$. 
In particular, our consideration is also a 
natural generalization of Bohr's theorem (Theorem 1).
Notice that it is shown in \cite{AM} that any power 
series (3), converging in $D$, converges in the sense 
of the left part of (7) for all $r$, $0<r<1$.
\begin{teo} 
The inequality 
$$1 - \sqrt[n]{\frac{2}{3}} < B_n (D). \eqno{(8)}$$
is true for any complete bounded $n$-circular domain $D$.
\end{teo}

This estimate can be improved for concrete domains.
\begin{teo}
For the unit ball
$D^1 = \{z: |z_1|^2 + \ldots + |z_n|^2 < 1 \}$ 
the following estimate is true:
$$B_n (D^1) >
\left\{\begin{array}{lll}
\frac{2}{5n} & \mbox{for} \quad n > 1 ,\\
&  \\
\frac{1}{2n} & \mbox{for large enough} \;  n.\\
\end{array}
\right.
$$
\end{teo}
\begin{teo}
For the unit hypercone $D^0 = \{z: |z_1| + \ldots +|z_n| < 1 \}$
the following inequality holds:
$$B_n (D^\circ ) < \frac{0,446663}{n} \; .
\eqno{(9)}$$
\end{teo}

\paragraph{Corollary.} 
$$1 - \sqrt[n]{\frac{2}{3}} < B_n ( D^\circ ) < 
\frac{0.446663}{n} \;. \eqno{(10)}$$

\paragraph{Remark 1.} 
The following asymptotic equality
$$1 - \sqrt[n]{\frac{2}{3}}= \frac{\log 3/2}{n} +
O (\frac{1}{n^2}) $$
is true, where $\log{3/2} \approx 0.405465$.  
Denoting by $B_{-}$ the left part of (10) and by $B_{+}$ 
the right one, we get
$$1 < \lim_{n \rightarrow \infty} \sup \frac{B_+}{B_{-}} < 1,1016. 
$$

Consider now another problem. Denote by $L_n(D^\circ )$ 
the biggest number $r$ 
such that if the series (3) converges in the hypercone
$D^\circ $ and (4) holds in it, 
then 
$$
\sum_{\alpha } \| c_{\alpha } z^{\alpha} \|_{L_1 
(\partial D^{\circ}_r )} \;
< \; 1 ,\eqno{(11)}
$$
where the $L^1$-norm is considered with respect to the measure 
${\mu}_r$, which is  the image of the measure
$$d \mu = \frac{(n-1)!}{(2 \pi i)^{n}} 
d | z_1 | \wedge \cdots \wedge d | z_{n-1}
| \wedge \frac{d z_1}{z_1} \wedge \cdots \wedge \frac{d z_n}{z_n}$$
by the homothetic transformation $z \to rz$. 
Usually the measure $d\mu$ is used
for calculating of Szego kernel for $D^0$ (see \cite{AY});
notice that $\mu(\partial \, D^0) = 1$. For $n=1$
this problem coincides with Bohr's problem. The analogous
problem for the polydisk $U_1$ (with $L^1$ on its Silov
boundary with respect to usual Lebesgue measure on it)
is equivalent to the problem, considered in Theorem~2.
\begin{teo}
For the hypercone $D^0$ the following estimates are true:
$$0,238843 < L_n ( D^\circ ) \leq 1/3 \; . \eqno{(12)}
$$\end{teo} 

It is surprising here that the estimates in (12) 
does not depend on $n$, which is different from the 
results of theorems 2---6.

Now consider a variant of multidimensional Bohr's problem,
dealing with expansions in series of homogeneous polynomials,
which is also a natural generalization power series expansion.
Let $Q$ be a complete circular domain (Cartan's domain) centered
at $0 \in Q$;  then any function $f(z)$, holomorphic in $Q$ can be 
expanded into the series
$$f (z) = \sum_{k=0}^\infty P_k (z),
\eqno{(13)}$$
where $P_k(z)$ is a homogeneous polynomial of degree $k$.
\begin{teo}
If the series (13) converges in the domain $Q$ and the estimate
$|f(z)| < 1$ holds in it, then
$$\sum_{k=0}^\infty | P_k (z) | < 1 \eqno{(14)}$$
in the homothetic domain $(1/3)  Q $.
Moreover, if $Q$ is convex,  then $1/3 $
is the best possible constant. 
\end{teo}

\section{Proofs}
\paragraph{Proof of Theorem 4.} The following generalization of Cauchy
inequalities was considered  in \cite{AM}: 
if (4) holds in $D$ then
$$| c_\alpha | < \frac{1}{d_\alpha (D)},\eqno{(15)}$$
where $d_{\alpha}(D) = \max\limits_{\bar{D}} {|z^{\alpha}|}$. 
Using Wiener's method it is easy to strengthen the estimates (15):
$$| c_\alpha | < (1 - | c_0 |^2)
\frac{1}{d_\alpha (D)} \eqno{(16)}$$
for $|\alpha| = {\alpha}_1 + \ldots + {\alpha}_n > 1$. 
We do not present the proof of (16), 
because it repeats the proof for polydisk $U_1$ 
(see \cite{BK}) but deals with inequalities (15)
instead of Cauchy inequalities. 
If (4) holds in $D$ then, applying (16), we get
\begin{eqnarray*}& & \sum_\alpha \sup_{D_r} | c_\alpha z^\alpha | =
\sum_\alpha | c_\alpha | d_\alpha (D_r) <
| c_0 | + \\ & + & (1 - | c_0 |^2 ) 
\sum_{|  \alpha | =1}^\infty
\frac{d_\alpha (D_r)}{d_\alpha (D)} =
| c_0 | + (1- | c_0 |^2)
\sum_{k=1}^\infty \left(\begin{array}{c} n+k-1\\ k \end{array} \right) 
r^k =\\ & = & | c_0 | + (1- | c_0 |^2) \biggl[ 
\frac{1}{(1-r)^n} -1 \biggr].\end{eqnarray*}

If now
$$\frac{1}{(1-r)^n} -1 \leq \frac{1}{2}, \eqno{(17)}$$
then
$$\sum_\alpha \sup_{D_r} | c_\alpha z^\alpha |
< | c_0 | + (1 - | c_0 |^2) \frac{1}{2} = 1 -
\frac{1}{2} (1 - | c_0 |)^2 < 1.$$
The condition (17) means that (7) is true if 
$r \leq 1-\sqrt[n]{\frac{2}{3}}$.

\paragraph{Proof of Theorem 5.} 
Consider 
Borel probability measure on $\partial D^1$, which is invariant under all
unitary transformations of $\C^n$: 
$$d \mu = \frac{(n-1)!}{(2 \pi i)^n} d | z_1 |^2 
\wedge \cdots \wedge  d | z_{n-1} |^2 \wedge
\frac{d z_1}{z_1} \wedge \cdots \wedge  \frac{d z_n}{z_n}.$$
The monomials $z^{\alpha}$ 
are orthogonal  with respect to integration by  $ \mu $
and 
$$\int_{\partial D^1} | z^{2 \alpha} | d \mu =
\frac{\alpha_1 ! \cdots \alpha_{n}! (n-1)!}{(| \alpha | + n-1)!} \;.
\eqno{(18)}$$

Further we repeat the proof of Theorem 2 
from \cite{BK} by Wiener method, but integrating on the sphere 
$\partial D^1$ by the measure $\mu$ instead
of integrating on the unit torus as  
in \cite{BK}. We obtain
$$\sum_{| \alpha |=k} | c_\alpha |^2
\frac{\alpha_n! \cdots \alpha_n! (n-1)!}{(| \alpha | + n-1)!}
< (1- | c_0 |^2)^2.$$

Furthermore, notice that from Schwarz Lemma 
for the ball $D^1$ it follows (again by using
Wiener's method from \cite{BK})
$$\sum_{| \alpha | =1} | c_\alpha |^2
< (1 - | c_0 |^2)^2.$$
Then, remembering that
$$d_\alpha (D^1) = \sqrt{\frac{\alpha_1^{\alpha_1} \cdots 
\alpha_n^{\alpha_n}}{| \alpha |^{| \alpha |}}},$$
where $0^0 = 1$, we get
\begin{eqnarray*}
&&  \sum_{\alpha} \sup_{D_r} | c_\alpha z^\alpha | =
| c_0 | + (\sum_{| \alpha | =1} | c_\alpha |) r +\\
& + & \sum_{k=2}^\infty \sum_{| \alpha | =k}
| c_\alpha 
| \sqrt{\frac{\alpha_1^{\alpha_1} \cdots
\alpha_n^{\alpha_n}}{| \alpha |^{| \alpha |}}}r^k 
< \\ & < & | c_0 | + (1- | c_0 |^2 ) r \sqrt{n}
+ (1 - | c_0 |^2) 
\sum_{k=2}^\infty
\biggl(\sum_{| \alpha | = k}
\frac{(k+n-1)! \alpha_1^{\alpha_1} \cdots \alpha_n^{\alpha_n}}
{k^k \alpha_1! \cdots \alpha_n! (n-1)!} \biggr)^{1/2}
r^k 
< \\
& < & | c_0 | + (1- | c_0 |^2) r \sqrt{n} +
(1- | c_0 |^2) \sum_{k=2}^\infty
\biggl( \frac{(k+n-1)!}{k! (n-1)!} \sum_{| \alpha | =k}
\frac{k!}{\alpha_1! \cdots \alpha_u!} \biggr)^{1/2} r^k = \\
& = & | c_0 | + (1- | c_0 |^2) r \sqrt{n} +
(1- | c_0 |^2)
\sum_{k=2}^\infty
\sqrt{\left(
\begin{array}{c} k+n-1\\ k \end{array} \right)} 
(r \sqrt{n})^k \end{eqnarray*}

It remains now to apply the estimate, obtained in Remark 1 from
\cite{BK}.

\paragraph{Proof of Theorem 6.} 
Consider the function
$$f_a (z) = \frac{1+a}{2} 
\frac{1- (z_1 + \cdots + z_n)}{1-a
(z_1 + \cdots + z_n)} =
\sum_\alpha c_\alpha z^\alpha,$$
where $0 < a < 1 $;  then $|f_a(z)| < 1$ in hypercone $D^\circ $ and 
$$\sum_\alpha | c_\alpha z^\alpha | =
\frac{1+a}{2} + \frac{1-a^2}{2}
\sum_{k=1}^\infty
a^{k-1} \sum_{| \alpha | =k}
\frac{k!}{\alpha_1! \cdots \alpha_n!} 
| z^\alpha |.
$$ Since
$$d_\alpha (D^\circ ) =\frac{\alpha_1^{\alpha_1} 
\cdots \alpha_n^{\alpha_n}}{|
\alpha |^{| \alpha |}},$$
it follows 
\begin{eqnarray*}
& & \sum_\alpha | c_\alpha | d_\alpha (D^\circ_r) =
\frac{1+a}{2} + \frac{1-a^2}{2}
\sum_{k=1}^\infty
\sum_{| \alpha | = k}
a^{k-1} \frac{k! \alpha_1^{\alpha_1} \cdots 
\alpha_n^{\alpha_n}}{\alpha_1!
\cdots \alpha_n! k^k} r^k>\\ & \geq  &
\frac{1+a}{2} + \frac{1-a^2}{2}
\sum_{k=1}^\infty a^{k-1} 
\frac{k!}{k^k} \biggl( \sum_{| \alpha | = k}
\frac{1}{\alpha_1! \cdots \alpha_n!} \biggr) r^k =\\ & = & 
\frac{1+a}{2}+ \frac{1-a^2}{2a} \sum_{k=1}^\infty
\frac{(a n r)^k}{k^k} \geq 1,\end{eqnarray*}
if
$$\sum_{k=1}^\infty \frac{(anr)^k}{k^k} \geq \frac{a}{1+a}.
$$

Let $x_0 (a)$ be the root of the equation
$$\sum_{k=1}^\infty \frac{x^k}{k^k} = \frac{a}{1+a};
$$
then, if $anr \geq x_0 (a)$, 
(7) fails for that $r$ and $D = D^0$. 
Now considering  
$a \to 1$, we obtain that (7) is not true for $D =D_0$ if 
$r \geq \frac{x_0}{n}$, where $x_0 = x_0(1)$.
Take notice that $x_0$ is a root of the equation
$$\sum_{k=1}^\infty \frac{x^k}{k^k} = 1/2, \eqno{(19)}$$
hence $B_n(D^0) \leq x_0 /n$. 
By using the program ``Mathematica 3.0'' \cite{W}, 
we estimated from above $x_0$. 
We obtained that the equation
$$\sum_{k=1}^p \frac{x^k}{k^k} = 1/2
\eqno{(20)}$$
has a root 0,446662 (where the last decimal digit is precise) 
if $p$ runs from 5 till 25. 
So the equation (20) has a root
less than 0,446663 if $5 \leq p \leq 25$, hence this estimate 
is true for $x_0$.

\paragraph{Proof of Theorem 7.} 
Take notice that by analogy with (18) it is true that
$$\int_{\partial D_r^\circ} | z^\alpha | d \mu_r =  
\frac{\alpha_1! \cdots \alpha_n! (n-1)!}{(| 
\alpha | + n-1)!} r^{| \alpha |} .$$

Using (16), we get
\begin{eqnarray*} & & \sum_\alpha \|
c_\alpha z^\alpha \|_{L^1 (\partial D^\circ_r)} <
| c_0 | + (1- | c_0 |^2) \sum_{k=1}^\infty
\sum_{| \alpha | = k}\\ & & \frac{k^k (n-1)!}{(n+k-1)!} 
\frac{\alpha_1! \cdots
\alpha_n!}{\alpha_1^{\alpha_1} \cdots
\alpha_n^{\alpha_n}}
r^{| \alpha |}
< | c_0 | + (1- | c_0 |^2 ) \times\\ & & \times 
\sum_{k=1}^\infty \frac{k^k (n-1)!}{(n+k-1)!}
\sum_{| \alpha | =k} r^{| \alpha |} =
| c_0 | + (1- | c_0 |^2) 
\sum_{k=1}^\infty \frac{k^k}{k!} r^k.\end{eqnarray*}

Denote by $r_0$ the root of the equation
$$\sum_{k=1}^\infty \frac{k^k}{k!} r^k = \frac{1}{2},
\eqno{(21)}$$
then (11) holds for $r = r_0$, therefore $L_n(D^0) \geq r_0$.
The left side of (21) is less than the left side
of the equation
$$\sum_{k=1}^{25} \frac{k^k}{k!} r^k + \frac{1}{\sqrt{52 \pi}} 
\frac{(ex)^{26}}{1-ex} =1/2,
\eqno{(22)}$$
what follows from an approximation of Gamma-function. 
Using again ``Mathematica 3.0'', we show that the root 
of the equation (22) is bigger than 0,238843, therefore 
the left part of the inequality (12) is true.

Further, consider, as in proof of Theorem 6, a function $f_a(z)$, 
for which
\begin{eqnarray*} & & \sum_\alpha \| c_\alpha z^\alpha 
\|_{L^1 (\partial D_r^\circ)} = \frac{1+a}{2} + 
\frac{1-a^2}{2} \sum_{k=1}^\infty 
\frac{k! (n-1)! a^{k-1}}{(n+k-1)!} \times\\
& & \times \sum_{| \alpha | = k}
r^{| \alpha |} =
\frac{1+a}{2} + \frac{1-a^2}{2a}
\sum_{k=1}^\infty (ar)^k \geq 1,
\end{eqnarray*}
if
$$\sum_{k=1}^\infty (ar)^k \geq \frac{a}{1+a}.$$
When  $a \to 1$, we get that the inequality (11) fails 
if $r > 1/3$.

\paragraph{Proof of Theorem 8.} 
In each section of the domain $Q$ by a complex line 
$$\alpha = \{z: z_j = a_j t, \quad j=1, \cdots, n; 
\quad t \in \C \}$$
the series turns into the power series by $t$
$$f (at) = \sum_{k=0}^\infty P_k (a) t^k$$
and, in addition, $|f(at)| < 1$. 
By Theorem 1
$$\sum_{k=0}^\infty | P_k (a) t^k | < 1
$$
in the section $\alpha \cap (\frac{1}{3} \cdot Q)$. 
But it is just (14), since 
$\alpha$ is an arbitrary complex line passing through
the origin.
Inversely, let the domain $Q$ be convex, then $Q$ is 
an intersection of half-spaces 
$$Q = \bigcap_{a \in J} \{ z: Re (a_1 z_1 + \cdots + a_n z_n)< 1 \}$$
with some $J.$ 
Since $Q$ is circular, we obtain
$$Q = \bigcap_{a \in J} \{z: | a_1 z_1 + \cdots + a_n z_n 
| < 1\} .$$

It is sufficient now to show that the constant $1/3$ cannot be
improved for each domain $P_a = \{z: |a_1z_1 + \ldots + a_nz_n |
< 1 \}$.  From Theorem 1 it follows that for any $r > 1/3$ there
exists a function $f(z)$, represented by (1) and such that
$|f(z)| < 1$ in the unit disk, but (2) fails in the disk $\{z:
|z| < r\}$.  To finish the proof we use the functions $f(a_1z_1 +
\ldots + a_nz_n)$.

\section{Final Remarks}
\paragraph{Remark 2.} In the proof of Theorem 5 the fact, that 
the domain in consideration is the ball 
$D^1 $ is used twice as follows: 
in order to prove that the monomials $z^\alpha $ are orthogonal 
on the sphere $\partial D^1 $ with respect to the measure 
$ d\mu $ and to get the estimate $d_{\alpha } (D^1 )  \leq 1 $.  
Therefore an analogous theorem holds for any 
complete $n$-circular domain $D$ if  $D^1 \subset D \subset U_1$.
But, for the hypercone $D^\circ $  this is not true 
(see Theorem 6).

\paragraph{Remark 3.}  In the proof of Theorem 6
was used the following quite rough inequality:
$$\sum_{| \alpha |=k}
\frac{\alpha_1^{\alpha_1} 
\cdots \alpha_n^{\alpha_n}}{\alpha_1!
\cdots \alpha_n!} \geq \sum_{| \alpha | = k}
\frac{1}{\alpha_1! \cdots \alpha_n!} .$$
Therefore, in fact, the estimate from above 
in the theorem might be decreased. 
For example, if $n = 2$ (10) implies $B_2(D^0) < 0.223332$, 
but using PC one can show that $B_2(D^0) < 0,191373$.
Hence,   $ 0.183502 < B_2 (D^\circ ) < 0.191373 $.

\paragraph{Remark 4.} Comparing Theorem 2 and Corollary we can see that
$B_n(D)$ depends essentially on the domain $D$ since $B_n(U_1) =
K_n$.  It seems that it is more natural to consider not a single
number in the problem, considered in Theorem 2 and Theorem 3, but
the largest subdomain $D_B$ of $D$ such that (5) holds.  From
\cite{BK} it follows, for example, that ${U_1}_B$ contains the
ball $\frac{1}{3} \cdot D^1$.

\paragraph{Remark 5.} Notice that
there exist such unbounded $n$-circular domains $D$ for which the
problems defined by the conditions (5) and (7), are equivalent.
Consider, for example, the domains $D = \{z: |z|^{\beta} < c\}$,
where $\beta$ is a multi-index with coprime components.  Bounded
holomorphic functions in such domains depend only on one variable
$z^{\beta}$ and the exact value of Bohr's radius equals to
$(\frac{1}{3})^{1/{|\beta |}}$.  Thus there exists $n$-circular
domains with Bohr's radius arbitrarily close to 1. Therefore it
is impossible to remove the assumption about convexity in
Theorem~8.

\paragraph{Remark 6.} Unlike Theorems 1---7, in Theorem 8
the series (13) is not a basis expansion. In \cite{A} it was
shown that there exists basis in the space of all holomorphic
functions in $Q$, consisting of homogeneous polynomials
$P_{k,m}(z)$, where $k$ is the degree of the polynomial and $m =
1, \ldots, \biggl( \begin{array}{c} k+n-1\\k\end{array} \biggr)$.
It is reasonable to consider Bohr's problem for such basis
expansions, but there is no results yet. This question is a
particular case of a more general problem in Bohr spirit for
expansions by an arbitrary basis in a domain $D \subset {\C^n} $
(under some restrictions on the basis, because, for example, the
basis $\{ 1$, $(z-1)/2$, $z^2$, $z^3$, $\ldots \} $ in the unit
disk has no Bohr's constant).

\section*{Acknowledgment} 
I am sincerely thankful to D.~Khavinson, who turned my attention
to Bohr's result (Theorem~1) and its first multidimensional
analogues (Theorems 2 and~3), and to V.~Zahariuta and P.~Djakov
for their help in preparing the English version of the paper.

The results of this paper were obtained during the author's 
visit in the Middle East Technical University.
The author would like to thank A.~Aytuna  and S.~Alpay for the 
invitation, and as well, A.~Aytuna and B.~Karas\"ozen for their 
support.

\hspace{15mm}

\begin{tabbing}
\parbox[t]{3in} {Address:\\
Department of Mathematics \\
and Computer Science, \\
Bar-Ilan University\\
52900 Ramat-Gan, ISRAEL\\
e-mail: aizenbrg@macs.biu.ac.il}

\parbox[t]{3in} {Temporary address:\\
Department of Mathematics\\
Middle East Technical University\\
06531 Ankara, TURKEY\\
e-mail: aizenbrg@arf.math.metu.edu.tr}
\end{tabbing}

\end{document}